
\magnification=\magstep1
\input amstex
\documentstyle{amsppt}
\TagsOnRight
\NoRunningHeads

\newdimen\notespace \notespace=7mm
\newdimen\maxnote  \maxnote=13mm

\define\marg#1 {\strut\vadjust{\kern-\dp\strutbox
  \vtop to\dp\strutbox{\vss \baselineskip=\dp\strutbox
  \moveleft\notespace\llap{\hbox to \maxnote{\hfil\rm#1}}\null}}}

\define\obs{\strut\vadjust{\kern-\dp\strutbox
  \vtop to\dp\strutbox{\vss \baselineskip=\dp\strutbox
  \moveleft\notespace\llap{\hbox to \maxnote{\hfil$!$}}\null}}\ }

\def\colim{\mathop{\vtop{\ialign{##\crcr
 \hfil\rm colim\hfil\crcr\noalign{\nointerlineskip}\rightarrowfill\crcr
 \noalign{\nointerlineskip\kern-.2326ex}\crcr}}}}

\define\Spec{\operatorname{Spec}}
\define\id{\operatorname{id}}

\define\A {\bold A}

\topmatter
\title 
Infinite Intersections of open subschemes and the Hilbert scheme of points.
\endtitle 
\author Roy M. Skjelnes and Charles Walter \endauthor
\affil 
Department of Mathematics, KTH and Laboratoire J.A. Dieudonn\'e, Nice
\endaffil
\address KTH, Stockholm, S-100 44, Sweden\endaddress
\email skjelnes\@math.kth.se \endemail 
\address 
Laboratoire J.-A.\ Dieudonn\'e (UMR 6621 du CNRS), Universit\'e de
Nice -- Sophia Antipolis, 06108 Nice Cedex 02, France \endaddress
\email walter\@math.unice.fr \endemail 
\abstract 
We study infinite intersections of open subschemes and the
corresponding infinite intersection of Hilbert schemes. If
$\{U_{\alpha} \}$ is the collection of open subschemes of a variety
$X$ containing the fixed point $P$, then we show that the Hilbert
functor of flat and finite families of $\Spec ({\Cal
O}_{X,P})=\bigcap_{\alpha} U_{\alpha}$ is given by the infinite
intersection $\bigcap _{\alpha}{\Cal Hilb}_{U_{\alpha}}$, where ${\Cal
Hilb}_{U_{\alpha}}$ is the Hilbert functor of flat and finite
families on $U_{\alpha}$. In particular we show that the Hilbert
functor of flat and finite families on $\Spec ({\Cal O}_{X,P})$ is
representable by a scheme.
\endabstract

\thanks The first author obtained financial support from Hellmuth
Hertz Foundation and  from Kungl. Vetenskapsakademien. \endthanks

\endtopmatter
\document

\head 1. - Introduction \endhead

We will consider in this article  infinite intersections of open
subschemes $\{U_{\alpha}\}$ of a fixed ambient scheme $X$. We are
interested in the corresponding Hilbert scheme and in particular in
the Hilbert scheme of $\Spec({\Cal O}_{X,P})$ the intersection  of the
open subschemes containing a point $P$ in $X$.

For a scheme $X$ the Hilbert scheme of $n$-points $\text{Hilb}^n_X$
(if it exists) represents the functor of finite flat families of
length-$n$ closed subschemes of $X$. Grothendieck constructed
$\text{Hilb}^n_X$ for $X$ quasi-projective over a noetherian base
scheme, but we wish to look at $\Spec ({\Cal O}_{X,P})$ for $P$ a
point of such an $X$. We know that if $U$ is an open subscheme of $X$
then $\text{Hilb}^n_U$ is an open subscheme of $\text{Hilb}^n_X$, so
there is a natural candidate for the Hilbert scheme of points on an
infinite intersection $\bigcap U_{\alpha}$ of open subschemes of $X$,
namely the corresponding infinite intersection $\bigcap
\text{Hilb}^n_{U_{\alpha}}$. Note (see Proposition (2.3)) however that
an infinite intersection of open subschemes is not necessarily a
scheme!

We restrict ourselves to infinite intersection of locally principal
open subschemes. The technical heart of the paper is the study of such
infinite intersections, which we call localized schemes. The notion of
localized schemes and generalized fraction rings is carried out in
Section (3). These concepts are thereafter applied to show that the
Hilbert functor of points on  a localized scheme ${\Cal S}^{-1}X$ is
representable, if  the Hilbert scheme of points on $X$ exists. A
special case of that statement gives the following.

\proclaim{Result} 
Let $X\to S$ be a projective morphism of noetherian schemes. Let $P$
be a point on $X$, with stalk ${\Cal O}_{X,P}$. Then the Hilbert
functor of $n$-points on $\Spec ({\Cal O}_{X,P})$ is representable by
a noetherian scheme $\operatorname{Hilb}^n_{\Cal
O_{X,P}}$. Furthermore, if $\{U_{\alpha }\}_{\alpha \in \Cal A}$ is
the collection of open subschemes of $X$ containing the point $P$,
then the Hilbert scheme of $n$-points on $\Spec({\Cal O}_{X,P})$ is
given as the infinite intersection
$$ \operatorname{Hilb}^n_{\Cal O_{X,P}}=\bigcap_{\alpha \in \Cal
A}\operatorname{Hilb}^n_{U_{\alpha}}.$$
\endproclaim

The above localization property for Hilbert functors of points was
known to hold for the affine line $X=\Spec(k[x])$ (see \cite{LS} and
\cite{S}) where the Hilbert scheme of points on fraction rings of
$k[x]$ were constructed explicitly. Here we show that the localization
property of the Hilbert functors of points hold for localized schemes.

What happens is the following. If $L$ is a line bundle on $X$, then we
get by pulling back $L$ to the universal family of $n$-points on $X$,
a vector bundle of rank $n$ over the Hilbert scheme
$\operatorname{Hilb}^n_X$. From each global section of $L$ we get a
determinant section of the norm bundle $N(L)$ on
$\operatorname{Hilb}^n_X$. If $U_s \subseteq X$ is the open subscheme
defined by the non-vanishing of a section $s\in \Gamma(X,L)$, then we
show that the Hilbert scheme of $n$-points on $U_s$ is the open
subscheme of $\operatorname{Hilb}^n_X$ given by the non-vanishing of
the corresponding determinant section of the norm bundle $N(L)$ on
$\operatorname{Hilb}^nX$.

\head 2. - Infinite Intersections of open subschemes \endhead

Let $X$ be a scheme, and let $\{U_\alpha \subseteq X \mid \alpha \in
{\Cal A}\}$ be a collection of open subschemes of $X$. The
set-theoretic intersection $\bigcap_{\alpha \in {\Cal A}} U_\alpha$
can be made into a locally ringed space by giving it the topology
induced by the Zariski topology of $X$ and by using as structural
sheaf the inverse image sheaf $i^{-1}\Cal O_X$, where $i : \bigcap
_{\alpha \in \Cal A}U_{\alpha} \to X$ is the inclusion.

In the category of locally ringed spaces we have that $\bigcap_{\alpha
\in {\Cal A}}U_{\alpha}=\varprojlim_{\alpha \in \Cal A}U_\alpha$. When
$\bigcap_{\alpha \in {\Cal A}} U_\alpha$ is a scheme, we also have
that $\bigcap_{\alpha \in \Cal A} U_\alpha = \varprojlim_{\alpha \in
A} U_\alpha$ in the category of schemes.  However, $\bigcap_{\alpha
\in \Cal A} U_\alpha$ is not necessarily a scheme; indeed
$\varprojlim_{\alpha \in A} U_\alpha$ does not always exist in the
category of schemes (see Proposition 2.3 below).

\proclaim{2.1. Theorem {\rm (Grothendieck)}} Let $\{
U_{\alpha}\}_{\alpha \in \Cal A}$ be a collection of open subschemes
of a scheme $X$. If the inclusion maps $i_\alpha : U_\alpha \to X$ are affine
morphisms, then the locally ringed space $\bigcap_\alpha
U_\alpha$ is a scheme.  Moreover, the inclusion $\bigcap_\alpha
U_\alpha \to X$ is an affine monomorphism.
\endproclaim

\demo{Proof}
All the assertions of the theorem are proven in \cite{EGA} {\bf
IV} \S 8.2 when the system $\{ U_\alpha
\}$ of open subsets is {\it filtered}, i.e.\ for any $\alpha$, $\beta$
there exists a $\gamma$ such that $U_\gamma \subseteq U_\alpha \cap
U_\beta$.  But we may reduce to the filtered case by replacing $\{
U_\alpha \}$ with the system of all finite intersections $\{
U_{\alpha_1} \cap \cdots \cap U_{\alpha_r} \}$ because the inclusion
maps remain affine morphisms while the categorical limit is
unchanged.

The construction of \cite{EGA} is that if the inclusions
$U_\alpha \subseteq X$ come locally from maps of commutative rings
$A \to B_\alpha$, then $\bigcap_\alpha U_\alpha \to X$ comes
from $A \to \colim_\alpha B_\alpha$.  We will use this in later
arguments.
\qed
\enddemo

\subsubhead 2.1.1. 
Locally principal subschemes \endsubsubhead An open subscheme
$U\subseteq X$ is {\it locally principal} if $X$ can be covered by
affine open subschemes $\Spec(A_i)$ such that each $U\cap \Spec(A_i)$
is a principal affine open subscheme of $\Spec(A_i)$ (i.e. of the form
$\Spec(A_{i,f_i})$ for some $f_i \in A_i$). The inclusion $U\subseteq
X$ of a locally principal open subscheme is an affine morphism, so
Theorem (2.1) applies.

\proclaim{2.2. Corollary} 
If the $U_\alpha \subseteq X$ are locally principal open subschemes
for all $\alpha \in \Cal A$, then the locally ringed space
$\bigcap_{\alpha \in \Cal A}U_{\alpha}$ is a scheme.
\endproclaim

\proclaim{2.3. Proposition} 
Let $X=\Spec(k[x,y])$ be the affine plane over a field $k$, and
consider the infinite intersection $\bigcap _{\alpha }(X\setminus
E_{\alpha})$, where $E_{\alpha}$ is any finite collection of closed
points in $X$. The locally ringed space
$$ \bigcap_{\alpha}(X\setminus E_{\alpha})
=\varprojlim_{\alpha}(X\setminus E_{\alpha})$$ 
is not a scheme.
\endproclaim

\demo{Proof} First we make some observations about the topology on
$Y=\bigcap_{\alpha}(X\setminus E_{\alpha})$. As a set $Y$ clearly is
the union $\{\xi\}\cup X_1$, where $\{\xi \}$ is the generic point of
the plane and where $X_1$ is the set of generic points of irreducible
plane curves. To describe the open subsets of $Y$ we use the fact that
the topology of $Y$ is induced by the Zariski topology on the affine
plane $X$. A proper closed subset of the affine plane $X$ is the union
of a curve $\Spec(k[x,y]/f)$ and a finite set of closed points
$E$. Since $Y\cap E=\emptyset$ we get that the open subsets of $Y$ are
of the form $Y\cap U_f$, where $U_f=\Spec (k[x,y]_f)$.

If $Y$ were a scheme, and consequently covered by affine schemes,
there would exist $f\in k[x,y]$ such that $Y\cap U_f$ is non-empty and
affine. We will show that this is impossible. Thus, we assume that $Y$
is a scheme and we let $f\in k[x,y]$ be such that $Y\cap U_f=\Spec
(B)$. Let $\varphi : k[x,y]_f \to B$ the homomorphism of rings that
corresponds to the morphism of affine schemes $Y\cap U_f \to U_f$.

A morphism of schemes $Z \to Y\cap U_f$ is equivalent with a morphism
$Z \to U_f$ that factors via $U_f \setminus E_{\alpha}$ for all finite
set of closed points $E_{\alpha}$. Let $S\subseteq k[x,y]$ be the
multiplicatively closed set of non-zero polynomials in the variable
$x$. Then clearly $S\cap P$ is non-empty for any maximal ideal $P
\subseteq k[x,y]$, and consequently we have a morphism $\Spec
(k[x,y]_{f,S})\to Y\cap U_f$. Similarily, we let $T\subseteq k[x,y]$
be the set of polynomials in the variable $y$ and thus obtain the
following diagram
$$\CD
k[x,y]_f @>{\varphi}>> B @>{\sigma}>> k[x,y]_{f,S} \cr
 @. @V{\tau}VV @VV{\iota}V \cr
@. k[x,y]_{f,T} @>{\epsilon}>> k(x,y)\cr
\endCD
$$
where the four maps $\sigma \varphi, \tau \varphi, \iota $ and
$\epsilon$ all are localization maps associated to the different
multiplicative subsets. We have that $\epsilon \tau\varphi =\iota
\sigma \varphi$. Now $\varphi$ is an epimorphism in the category of
commutative rings because $Y\cap U_f \to U_f$ is clearly  a
monomorphism of affine schemes,  hence $\epsilon \tau =\iota \sigma$
and the diagram is commutative. It follows that the image of $B$ in
$k(x,y)$ lies in the intersection of the two subrings
$k[x,y]_{f,S}\cap k[x,y]_{f,T}$. We have that  $S\cap T=k^{*}$, the
set of non-zero constants of $k$, and we get that $k[x,y]_{f,S}\cap
k[x,y]_{f,T}=k[x,y]_f$. That is, the homomorphism $\varphi : k[x,y]_f
\to B$ has a retraction and dually the morphism of affine schemes
$Y\cap U_f \to U_f$ has a section. It follows that that the morphism
of schemes $Y\cap U_f \to U_f$ is surjective on points. The underlying
set of points of $Y=\{\xi\}\cup X_1$, where $X_1$ are the generic
points of the affine plane curves. But, since $U_f$ is assumed
non-empty the set $U_f$ contains closed points and consequently the
inclusion $Y\cap U_f \to U_f$ can not be surjective. We have reached a
contradiction and hence proven that no nonempty $Y\cap U_f$ is
affine. Since the open subsets of $Y$ are on the form $Y\cap U_f$ we
get that $Y$ contains no nonempty affine subschemes, and in particular
$Y$ is not a scheme. \qed
\enddemo

\subhead 2.4. Infinite intersection of Noetherian schemes \endsubhead

If $B$ is an $A$-algebra we denote with $IB$ the extension of an ideal
$I\subseteq A$, and with $J\cap A$ the contraction of an ideal
$J\subseteq B$.

\proclaim{2.5 Lemma} 
Let $\varphi : A\to B$ be a homomorphism of commutative rings. Assume
that the corresponding morphism of schemes $\Spec(B) \to \Spec (A)$ is
an open immersion. Then any ideal $J\subseteq B$ is the extension of
its contraction to $A$.
\endproclaim

\demo{Proof} 
Let $J\subseteq B$ be an ideal. The extension of the contraction
$(J\cap A)B$ is trivially contained in $J$ and we need only to show
that $J\subseteq (J\cap A)B$.

Since affine schemes are quasi-compact, $\Spec (B)$ can be covered by
a finite number of principal affine open subschemes of $\Spec
(A)$. Thus there exist $f_1, \ldots, f_r$ in $A$ such that $\Spec
(B)=\bigcup_{i=1}^r\Spec(A_{f_i})$.  The induced maps $A_{f_i} \to
B_{\varphi (f_i)}$ are isomorphisms, and one deduces that for any
element $x$ in the ideal $J\subseteq B$ there exist elements $a_1,
\ldots, a_r $ in $A$ such that $\varphi (a_i)=\varphi (f_i)^Nx$, for
some $N\gg 0$. In particular we have that $a_i \in \varphi
^{-1}(x)\subseteq J\cap A$ for each $i=1, \ldots, r$. Since the
$\Spec(B_{\varphi(f_i)})$ cover $\Spec (B)$ it follows that there
exist $b_1, \ldots , b_r$ such that $\sum_{i=1}^rb_i\varphi(f_i)^N
=1$. Then we have that $x=\sum_{i=1}^rb_i\varphi (a_i)$ is in the
extension of $J\cap A$, hence $J\subseteq (J\cap A)B$.
\qed
\enddemo

\proclaim{2.6. Lemma} 
Suppose we are given a direct (or filtered) system of commutative
rings $A_i$ and transition maps $\varphi_{ij} : A_j\to A_i$ such that
any ideal in $A_i$ is the extension of the contraction to $A_j$. Then
any ideal in $\colim_{i}A_i$ is the extension of the contraction to
$A_j$ via the natural homomorphism $A_j \to \colim_{i}A_i$.
\endproclaim
\demo{Proof} 
Let $J$ be an ideal of the direct colimit $A=\colim _{i}A_i$. From the
the assumption of the transition maps $\varphi_{ij}$ we have $ J\cap
A_i=(\varphi^{-1}_{ij}(J\cap A_i))A_i=(J\cap A_j)A_i$. By the
exactness of the direct colimit we get that $\colim_{i}(J\cap A_i)$ is
an ideal in $\colim_{i}A_i$, easily seen to coincide with $J$.\qed
\enddemo

\proclaim{2.7. Proposition} 
In the situation of Theorem (2.1), if $X$ is a noetherian scheme then
so is $\bigcap_{\alpha \in\Cal A}U_{\alpha}$.
\endproclaim

\demo{Proof} 
Assume first that $X=\Spec (A)$ is affine. Then $\bigcap _{\alpha \in
\Cal A}U_ {\alpha}\to X$ is given by $A\to
\colim_{\alpha}B_{\alpha}$. We must show that
$\colim_{\alpha}B_{\alpha}$ is noetherian. By Lemma (2.5) we have that
the homomorphism of rings $\varphi_{\alpha} : A \to B_{\alpha}$ is
such that the extension of the contraction of an ideal $J\subseteq
B_{\alpha}$ equals $J$. It follows from Lemma (2.6) that any ideal
$J\subseteq \colim_{\alpha}B_{\alpha}$ is the extension of its
contraction to $A$. Since $A$ is noetherian and consequently any ideal
of $A$ is finitely generated, it follows that any ideal of
$\colim_{\alpha}B_{\alpha}$ is finitely generated. Hence
$\colim_{\alpha}B_{\alpha}$ is noetherian.

If $X$ is simply a noetherian scheme, then $Y=\bigcap _{\alpha \in
\Cal A}U_{\alpha}$ is given locally by the construction above, so $Y$
is locally noetherian.  Since $X$ is quasi-compact and the morphism $Y
\to X$ is affine and hence quasi-compact, $Y$ is also quasi-compact.
Hence $Y$ is a noetherian scheme.\qed
\enddemo

\head 3. - Localized schemes and generalized fraction rings\endhead

\subhead 3.1. 
Localized schemes \endsubhead Let $X$ be a scheme. We will write
sections of invertible sheaves on $X$ as pairs $(s,L)$, where $s :
{\Cal O}_X \to L$ is a global section of the invertible sheaf $L$. We
let $U_s \subseteq X$ denote the open subscheme where the section $s$
is non-vanishing, that is the complement of the support of $s$.

\proclaim{3.2. Theorem} 
Let ${\Cal S}=\{(s_{\alpha},L_{\alpha})\}_{\alpha \in \Cal A}$ be a
collection of sections of invertible sheaves on $X$. Then there exists
a morphism of schemes $i_{\Cal S} : {\Cal S}^{-1}X \to X$ such that
the following two assertions hold.
\roster
\item 
The pull-back $i_{\Cal S}^{*}(s_{\alpha}) : {\Cal O}_{\Cal S ^{-1}X}
\to i_{\Cal S}^{*}L_{\alpha}$ is nowhere vanishing on ${\Cal
S}^{-1}X$, for all $\alpha \in \Cal A$.
\item  
Any homomorphism $f: T\to X$ of schemes such that $f^{*}(s_{\alpha}) :
{\Cal O_T} \to f^*L_{\alpha}$ is nowhere vanishing on $T$ for all
$\alpha \in \Cal A$, has a unique factorization via $i_{\Cal S}$.
\endroster
Moreover, $i_{\Cal S} :{\Cal S}^{-1}X \to X$ is unique up to unique
isomorphism.
\endproclaim
\demo{Proof} 
Each $U_{s_\alpha}\subseteq X$ is a locally principal open subscheme,
thus by Corollary (2.2) we have that the inclusion $\bigcap_{\alpha
\in \Cal A}U_{s_{\alpha}}\to X$ is a morphism of schemes. Let ${\Cal
S}^{-1}X=\bigcap _{\alpha \in \Cal A}U_{s_{\alpha}}$ and let $i_{\Cal
S}$ be the inclusion ${\Cal S}^{-1}X\to X$.

Note that ${\Cal S}^{-1}X=\varprojlim_{\alpha \in \Cal A}U_{s_\alpha}$
such that a morphism of schemes $f : T \to X$ factors via $i_{\Cal S}
: {\Cal S}^{-1}X \to X$ if and only if $f$ factors via $i_{s_{\alpha}}
: U_{s_{\alpha}} \to X$ for all $\alpha \in \Cal A$. Assertion (1)
then follows since it is clear that the pull-back of a section $s :
{\Cal O}_X \to L$ along the inclusion $i_s : U_s \to X$ is
non-vanishing.

To show Assertion (2) it suffices to show that for a given section
$(s,L)$ on $X$ a morphism $f : T\to X$ factors via $i_s : U_s \to X$
if and only if $f^{*}(s)$ is non-vanishing on $T$. We can cover $X$ by
open affine subschemes $\{\Spec (A_i)\}_{i\in \Cal I}$, such that $U_s
\cap \Spec (A_i)$ is given by some principal open subschemes $\Spec
(A_{i,f_i})$ of $\Spec (A_i)$. Assertion (2) now follows from the
universal properties of fraction rings.

It is clear that the condition on the morphism $f : T\to X$ given in
Assertion (2) defines a functor which is represented by the scheme
$\Cal S^{-1}X$ with universal element $i_{\Cal S} : {\Cal S}^{-1}X \to
X$, hence uniqueness follows.\qed
\enddemo

\proclaim{3.3. Lemma} 
Let ${\Cal S}=\{(s_{\alpha},L_{\alpha})\}_{\alpha \in \Cal A}$ be a
collection of sections of invertible sheaves on $X$, and $p^*{\Cal
S}=\{(p^{*}(s_{\alpha}),p^*L_{\alpha})\}_{\alpha \in \Cal A}$ the
pull-back of ${\Cal S}$ along a given morphism of schemes $p : Z\to
X$. Then the localization map $i_{p^*{\Cal S}} :(p^*{\Cal S})^{-1}Z
\to Z$ and the pull-back ${\Cal S}^{-1}X\times_XZ \to Z$ of the
localization map on $X$, coincide up to unique isomorphism.
\endproclaim
\demo{Proof} 
One immediately checks that the map ${\Cal S}^{-1}X\times_XZ \to Z$
satisfies the two conditions (1) and (2) of Theorem (3.1), which
proves the claim.\qed
\enddemo

\remark{3.3.1. Remark} 
Let $X$ be a scheme over some base $S$, and let $T \to S$ be a
morphism of schemes. Let ${\Cal S}=\{(s_{\alpha},L_{\alpha})\}_{\alpha
\in \Cal A}$ be a collection on $X$, and let $U_{s_\alpha}$ be the
locally principal open subscheme defined by the section
$(s_{\alpha},L_{\alpha})$. The natural morphism of schemes
$$ \bigcap_{\alpha \in \Cal A}(U_{s_{\alpha}}\times_ST)
=\varprojlim_{\alpha \in \Cal A}(U_{s_{\alpha}}\times_S T) \to
(\varprojlim_{\alpha \in \Cal A}U_{s_{\alpha}})\times_S T =(\bigcap
_{\alpha \in \Cal A}U_{\alpha})\times_ST
$$ 
is an isomorphism by Lemma (3.3).
\endremark

\subhead 3.4. 
Generalized fraction rings \endsubhead Let $R$ be a ring (commutative
with unit), and let $U=\{(s_{\alpha},L_{\alpha})\}_{\alpha \in \Cal
A}$ be a collection of pairs $s_{\alpha} \in L_{\alpha}$ with
$L_{\alpha}$ an invertible $R$-module. Let $\bold{N}\cdot \Cal A$
denote the subset of $\bold N^{\Cal A}$ consisting of systems of
non-negative integers $a=\{a_{\alpha}\}_{\alpha \in\Cal A}$ having
only a finite number of non-zero components. The set $\bold{N}\cdot
\Cal A$ is naturally partially ordered where we say that $a\leq b$ if
for each component we have $a_{\alpha} \leq b_{\alpha}$. We define for
any $a\in \bold{N}\cdot \Cal A$ the invertible $R$-modules
$$ 
\underline L^{a} =\bigotimes_{a_{\alpha}\neq 0}L_{\alpha}^{\otimes
a_{\alpha}} \quad \text{ and } \quad \underline L ^{0}=R.\tag{3.4.1}
$$ 
For any $b\in \bold{N}\cdot \Cal A$ we have a natural identification
$\underline L^{a} \otimes_R \underline L^{b} =\underline L^{a+b}$. We
have furthermore the element
$$
s^{b}=\bigotimes_{b_{\alpha} \neq 0}  s_{\alpha}^{\otimes
b_{\alpha}} \in \underline L^{b}\tag{3.4.2}
$$ 
The element $s^b$ defines $R$-module homomorphisms
$$ 
\underline L^{a}\to \underline L^{a+b} \tag{3.4.3}
$$
sending $x\in \underline L^{a}$ to $x\otimes s^{b} \in \underline
L^{a}\otimes _R \underline L^{b}=\underline L^{a+b}$. We denote the
direct colimit of the $R$-modules (3.4.1) and the described transition
maps (3.4.3) as
$$ 
R_U =\colim_{a\in \bold{N}\cdot \Cal A} \{\underline
L^{a}\}. \tag{3.4.4}
$$
Note that we have a natural product structure on $R_U$ as
$$ 
\underline L^{a} \cdot \underline L^b \subseteq \underline
L^{a+b},
$$ 
by $x^a\cdot y^b:=x^a\otimes y^b$. As the $R$-modules $\underline
L^{a}$ are invertible for all $a\in |\Cal A |$, we have that
$x^{a}\cdot y^{b}=y^b\cdot x^a$. Hence $R_U$ is a commutative
ring. As $R=\underline L^0$, we have that $R_U$ is a commutative
$R$-algebra. We call $R_U$ a generalized fraction ring (with respect
to $U=\{(s_{\alpha},L_{\alpha})\}_{\alpha \in \Cal A}$). If we have
$L_{\alpha}=R$ for all $\alpha$, the direct colimit $R_U$ is the
fraction ring $V^{-1}R$, where $V\subseteq R$ is the multiplicative
system generated by the $s_{\alpha}$.

\subhead 3.5. 
Properties of the generalized fraction rings \endsubhead Let
$U=\{(s_{\alpha},L_{\alpha})\}_{\alpha \in \Cal A}$ be a collection of
invertible modules. We will in this section list some properties of
the generalized fraction rings $R_U$, properties that we will use
later on in Section (3).

\remark{3.5.1. Remark} 
We have that $R_U$ is an $R$-algebra, thus also an $R$-module. By
definition $R_U$ is the direct colimit of locally free, in particular
flat, $R$-modules $\underline L^{a}$, hence $R_U$ is a flat
$R$-module.
\endremark

\remark{3.5.2. Remark} If $N$ is a $R$-module we denote by
$N_U:=R_U\otimes_R N$. We have that tensor product commute with direct
colimit  hence
$$ N_U = \colim_{a\in \bold{N}\cdot \Cal A}\{ \underline L^{a} \}\otimes_R N
=\colim_{a\in \bold{N}\cdot \Cal A}\{\underline L^{a}\otimes_R N\}.$$

In particular we have the following. Let $R$ be an $A$-algebra, and $A
\to B$ a homomorphism of rings. Write $R\otimes_AB=R_B$ and let $U_B$
be the collection on $R_B$ coming from the collection $U$ on $R$, that
is $U_B=\{(s_{\alpha}\otimes 1,L_{\alpha}\otimes_A 1)\}_{\alpha \in
\Cal A}$. Then we have that
$$ R_U\otimes_A B=\colim_{a\in \bold{N}\cdot \Cal A}(\underline
L^a)\otimes_A B=\colim_{a \in \bold{N}\cdot \Cal A}(\underline L^a
\otimes_A B)=(R_B)_{U_B}. \tag{3.5.2.1}$$
\endremark

\remark{3.5.3. Remark} 
Let $N$ be an $R$-module. For any element $x\in \underline
L^{a}\otimes_R N$ we denote the image of $x$ in the colimit $N_U$ by
$x/s^a$, where $s^a$ is the element defined in (3.4.2).  If $y\in
\underline L^b\otimes_R N$ is another element then $x/s^a=y/s^b$ in
$N_U$ if and only if there exists $c\in \bold{N}\cdot \Cal A$ such
that
$$ s^c (s^bx-s^ay)=0 \quad \text{ in } \underline L^{a+b+c}\otimes_R N.$$
In particular we have that $s_{\alpha}\in L_{\alpha}$ becomes a unit in
$R_U$, namely $s_{\alpha}/s_{\alpha}= 1$.
\endremark

\remark{3.5.4. Remark} 
An invertible $R$-module $\underline{L}^b$ is faithfully flat, hence a
map
$$ s^{a} : N \to \underline L^{a} \otimes_R N \tag{3.5.4.1}$$
is injective or surjective if and only if 
the $R$-module maps
$$ 
s^{a} : \underline L^b \otimes_R N \to \underline L^{a+b}\otimes_RN
\tag{3.5.4.2}
$$ 
is injective or surjective, respectively.
\endremark

\remark{3.5.5. Remark} For any subset $J\subseteq \Cal A$ we can
consider the colimit 
$$ R_{U_J}=\colim_{a\in \bold{N}\cdot J}\{\underline L^a\}.$$
The union of two subsets $J_1$ and $J_2$ of $\Cal A$ again is a
subset of $\Cal A$, and  we have that $\Cal A$ is partially ordered by
the union of
its subsets. It is clear that $R_U$ is the direct colimit
$$ 
R_U =\colim_{\text{ finite } J\subseteq \Cal A}
\{ R_{U_J} \}.\tag{3.5.5.1}
$$
\endremark
  
\proclaim{3.6. Proposition} 
Let $U=\{(s_{\alpha},L_{\alpha})\}_{\alpha \in \Cal A}$ be a
collection of invertible modules on $R$. We have the following.
\roster
\item For any $\alpha \in \Cal A$ we have that $1\otimes s_{\alpha}
\in R_U \otimes _R L_{\alpha}$ is nowhere vanishing.
\item If $R\to A$ is an $R$-algebra homomorphism such that $1\otimes
s_{\alpha}$ in $A\otimes _R L_{\alpha}$ is nowhere vanishing, for all
$\alpha \in \Cal A$, then the homomorphism $R \to A$ factors via the
homomorphism $R \to R_U$.
\endroster
\endproclaim

\demo{Proof} 
To show Assertion (1) we need to show that the map $R_U \to R_U
\otimes_R L_{\alpha}$ determined by sending $1$ to $1\otimes
s_{\alpha}$ is an isomorphism, for all $\alpha \in \Cal A$. We have
(3.5.2) that $R_U \otimes_R L_{\alpha}=R_U$, where we identify
$1\otimes s_{\alpha}$ with $s_{\alpha}$. We have already remarked
(3.5.3) that $s_{\alpha}$ is a unit in $R_U$ for all $\alpha \in \Cal
A$, and Assertion (1) follows.

We then show Assertion (2). From the assumption we have that $A \to
A\otimes_R L_{\alpha}$ sending $x\to x\otimes s_{\alpha}$ is an
isomorphism of $A$-modules for all $\alpha\in \Cal A$. It follows that
$A\to A\otimes \underline L^{a}$ is an isomorphism for all $a\in
\bold{N}\cdot \Cal A$, hence the colimit $R_U\otimes_RA$ is isomorphic
to $A$.  We have an $R$-algebra homomorphism $R_U \to R_U\otimes_R A$
that composed with the inverse of the isomorphism $A\to R_U
\otimes_RA$ gives our desired map. \qed
\enddemo

\proclaim{3.7. Corollary} 
Let $\Cal S=\{(s_{\alpha}, \tilde L_{\alpha})\}_{\alpha \in
\Cal A}$ be a collection of invertible sheaves on a affine scheme
$X=\Spec (R)$. Let $L_{\alpha}=\Gamma (X,\tilde L_{\alpha})$. Then $i
: \Cal S^{-1}X \to X$ is canonically identified with $\Spec (R_U) \to
X$.
\endproclaim

\demo{Proof} By the Proposition we have that $\Spec (R_U)\to \Spec
(R)$ satisfies the universal defining properties of ${\Cal S}^{-1}X
\to X$.\qed
\enddemo

\remark{3.7.1. Remark} If the collection $U=(s, L)$ consists of one
pair only then we have  
$ \Spec (R_{(s,L)})=U_s$, where $U_s\subseteq \Spec (R)$ is the
locally principal affine open subscheme defined by the
non-vanishing of the section $s\in L$.
\endremark

\remark{3.7.2. Remark} If the collection $U=\{(s_i,L_i)\}_{i=1,
\ldots, r}$ is finite then we
can reduce the situation to the single pair $(s,L)$, where
$$s=s_1\otimes \dots \otimes s_r \in L_1\otimes_R \cdots \otimes _R
L_r=L.$$
Then we have that $R_U=R_{(s,L)}$. On the level of Spec we have that
the finite intersection of locally principal open subschemes
$U_{s_i}\subseteq \Spec (R)$ is the locally principal open subscheme
$U_s\subseteq \Spec (R)$.
\endremark

\proclaim{3.8. Proposition} 
Let $N$ be an $R$-module, and $U$ a collection of invertible
$R$-modules. If the map $N \to N_U$ is an isomorphism of $R$-modules,
then the maps $N \to \underline L^{a}\otimes_A N$ are isomorphisms for
all $a \in \bold{N}\cdot \Cal A$.
\endproclaim

\demo{Proof} 
As $N_U$ is the direct colimit (3.5.2) it is clear that the assumed
injectivity of $N \to N_U$ implies that the maps $N \to \underline
L^{a}\otimes_R N$ are injective for all $a\in \bold{N}\cdot \Cal
A$. In particular the maps (3.5.4.2) are injective.  We need only to
show surjectivity of the maps $N \to \underline L^a\otimes_RN$. Let
$x\in \underline L^a\otimes_R N$. The map $N \to N_U$ to the direct
colimit is assumed to be surjective. Hence there exists $y\in N$
having the same image as $x$ in $N_U$. Thus $y=x/s^a$. By (3.5.3) we
have that there exists $c\in \bold{N}\cdot \Cal A$ such that
$$
s^c(ys^a-x) = 0 \quad \text{ in } \underline L^{a+c}\otimes_R N.
$$
As the maps (3.5.4.2) are injective we have that $ys^a=x$ in
$\underline L^{a}\otimes_RN$, hence we have proven the surjectivity of
$N \to \underline L^a\otimes_R N$. \qed
\enddemo

\proclaim{3.9. Lemma} 
Let $f: M \to N$ be a homomorphism of $R$-modules. Assume that $N$ is
finitely generated and that the induced map $N \to N_U$ is an
isomorphism of $R$-modules. If the $R_U$-linear map $f_U : M_U \to
N_U$ is surjective, then the homomorphism $f : M \to N$ is surjective.
\endproclaim

\demo{Proof} 
Let $x_1, \ldots ,x_r$ generate the $R$-module $N$. For each $a \in
\bold{N}\cdot \Cal A$ we let $f_a : \underline L^a \otimes_A M\to
\underline L^a\otimes_A N$ denote the induced $R$-linear maps. The map
$f_U$ between the direct colimits is assumed to be surjective, hence
there exists $a\in \bold{N}\cdot \Cal A$ and elements $y_1, \ldots
,y_r$ in $\underline L^a\otimes_R M$ such that $f_a (y_i)\in
\underline{L}^a\otimes_R N$ has the same image as $x_i$ in $N_U$ for
each $i=1, \ldots, r$. By Lemma (3.8) we have that $N \to \underline
L^a\otimes_R N$ is surjective, thus the images of $x_1, \ldots ,x_r$
generate $\underline L^a\otimes_RN$. It follows that the $R$-module
homomorphism $f_a : \underline L^a \otimes _RM\to \underline
L^a\otimes_RN$ is surjective. As $\underline L^a$ is a faithfully flat
$R$-module we obtain that $M \to N$ is surjective.\qed
\enddemo

\proclaim{3.10. Lemma} 
Let $U=\{(s_{\alpha}, L_{\alpha})\}_{\alpha \in \Cal A}$ be a
collection of invertible $R$-modules. Let $I_U\subseteq R_U$ be an
ideal of the generalized fraction ring $R_U$ and let $I=I_U\cap R$
denote its contraction. Then the localization map $R/I \to R_U/I_U$ is
an isomorphism if and only if $R_U/I_U$ is finitely generated as an
$R$-module.
\endproclaim

\demo{Proof} 
The only non-trivial part of the Lemma is to show that if $R_U/I_U$ is
finitely generated then $R/I \to R_U/I_U$ is surjective. We claim that
the ideal $I_U$ in $R_U$ is the extension of its contraction $I$.  By
(3.5.5.1) we have that $R_U$ is the direct colimit of generalized
fraction rings $R_{U_J}$, where $J$ is a finite subset of $\Cal A$. It
follows from (3.7.2) that $\Spec (R_{U_J}) \to \Spec (R)$ is an open
immersion. By Lemma (2.6) we have that any ideal in $R_{U_J}$ is the
extension of its contraction to $R$. It then follows by Lemma (2.7)
that any ideal $I_U$ in $R_U$ is the extension of its contraction to
$R$.

Thus we have that $R_U\otimes_R R/I=R_U/I_U$. Hence we may assume that
$R/I=R$. We must show that if $R_U$ is finitely generated $R$-module
then $R\to R_U$ is surjective. This is a special case of Lemma
(3.9).\qed
\enddemo

\proclaim{3.11. Proposition} 
Let $X\to S$ be a scheme over a base scheme $S$, and let ${\Cal S}$ be
a collection of sections of invertible sheaves on $X$. Let $f: T\to S$
be a morphism of schemes, and let $j: Z\subseteq {\Cal
S}^{-1}X\times_ST$ be a closed subscheme such that the projection map
$Z \to T$ is finite. Then $Z$ is a closed subscheme of $X\times_ST$
via the composite map $(i_{\Cal S}\times \id_T)\circ j$.
\endproclaim

\demo{Proof} 
We may assume that $T=\Spec (A)$ is affine since closedness is a local
property. We may also, by Lemma (3.3) assume that $T=S$. Finally, it
is clear that we may assume that $X=\Spec (R)$ is affine. The
Proposition now follows from Lemma (3.10).\qed
\enddemo

\head 4. - Determinants and Localized Schemes \endhead

There exists a notion of noncommutative localization and
$\sigma$-inverting rings, for any ring $R$ and any set $\sigma$ of
morphisms $s : P \to Q$ of finitely generated projective modules $P$
and $Q$ (\cite{C}, \cite{NR}). We will our commutative situation
obtain those $\sigma$-inverting rings as generalized fraction rings of
a collection of determinants and norm bundles.

\subhead 4.1. Notation \endsubhead
Let $s : E \to L$ be an $A$-module homomorphism between two locally
free $A$-modules $E$ and $L$ of finite rank $n$. We take the highest
exterior power of the $A$-module map $ s : E \to L$ and obtain en
element
$$ 
\text{det}(s)=\wedge s \in \text{Hom}_A(\wedge^n E, \wedge^n
L)=N(E,L).
$$
The element $\text{det}(s)$ is an element of the invertible $A$-module
$N(E,L)$. We clearly have that $s : E \to L$ is an isomorphism if and
only if $\text{det}(s): A \to N(E,L)$ is nowhere vanishing.

Let $\varphi : A \to B$ be an $A$-algebra homomorphism, and let
$E_B=E\otimes_AB$ and $s_B=s\otimes 1\in L_B=L\otimes_AB$. Then we
have that
$$
\text{det}(s_B)=\text{det}(s)\otimes_A1=\varphi (\text{det}(s)).
\tag{4.1.1}$$

Let now $E$ be an $A$-algebra such that $E$ is locally free of finite
rank as an $A$-module. Let $U=\{(s_{\alpha}, L_{\alpha})\}_{\alpha \in
\Cal A}$ be a collection of elements $s_{\alpha}$ in invertible
$E$-modules $L_{\alpha}$. We denote by
$$ 
N_{E}(U)=\{(\text{det}(s_{\alpha}), N(E,L_{\alpha}))\}_{\alpha \in
\Cal A} 
$$ 
the corresponding collection on $A$. If $U$ is a collection of
invertible modules on $E$ we refer to $N_{E}(U)$ as the corresponding
collection of norms on $A$.

\proclaim{4.2. Proposition} Let $E$ be an $A$-algebra such that $E$ is
locally free of finite rank $n$ as an $A$-module. Let $U$ be a
collection on $E$ and let $N_{E}(U)$ be the corresponding collection of
norms on $A$. For any homomorphism of rings $\varphi : A \to B$ the
following two statements are equivalent.
\roster
\item The induced homomorphism $E\otimes_A B \to E_U \otimes _A B$ is
an isomorphism
\item The homomorphism $\varphi : A \to B $ factors via $A\to
A_{N_{E}(U)}$.
\endroster
In particular we have that $E\otimes_A A_{N_{E}(U)} \to E_U\otimes_A
A_{N_{E}(U)}$ is an isomorphism.
\endproclaim

\demo{Proof} 
By (3.5.2.1) we have $E_U\otimes_A B=(E\otimes_AB)_{U_B}$, where $U_B$
is the collection $\{(s_{\alpha}\otimes 1,
L_{\alpha}\otimes_AB)\}_{\alpha \in \Cal A}$. The Assertion (1) then
reads by Propositon (3.7) that the sections $s_{\alpha}\otimes _A 1$
are nowhere vanishing, for all $\alpha \in \Cal A$. Hence their
determinants $\det(s_{\alpha}) \in N(E,L_{\alpha})$ are nowhere
vanishing. It then follows by the universal property of the
generalized fraction rings, Proposition (3.6), that the homomorphism
$f: A \to B$ factors via $A \to A_{N_E(U)}$. We have proven that
Assertion (1) implies Assertion (2). Assume now that Assertion (2)
holds. By Proposition (3.6) we have that the sections $1\otimes
\text{det}(s_{\alpha}) \in A_{N(U)}\otimes_A N(E,L_{\alpha})$ are
nowhere vanishing for all $\alpha \in \Cal A$. Then we have that
$f(1\otimes \text{det}(s_{\alpha})) \in B$ are invertible, for all
$\alpha \in \Cal A$. It follows from (4.1.1) that the sections
$s_{\alpha}\otimes 1$ in $L_{\alpha}\otimes_A B$ are nowhere
vanishing, for all $\alpha \in \Cal A$. Consequently $E_B=E\otimes_AB$
is isomorphic to the direct colimit $(E_B)_{U_B}$, which by (3.5.2.1)
equals $E_U\otimes_A B$.\qed
\enddemo

\definition {4.3. Definition} 
We say that a flat and finite morphism of schemes $q : Z \to H$ is of
relative rank $n$, if the quasi-coherent ${\Cal O}_H$-module $q_*{\Cal O}_Z$
is locally free of finite rank $n$.
\enddefinition

\subhead 4.4. 
Determinant sections \endsubhead Let $s : {\Cal O}_Z \to L$ be a
section of an invertible sheaf $L$ on $Z$. Let $q : Z \to H$ be a
morphism of schemes that is flat, finite and of relative rank $n$. We
then have that $q_*L$ is a quasi-coherent ${\Cal O}_H$-module, locally free
of rank $n$. The highest exterior power of the ${\Cal O}_H$-module
homomorphism $q_*(s) : q_*{\Cal O}_Z \to q_*L$ gives a global section
$\text{det}(s)$ of the invertible ${\Cal O}_H$-module
$$ \Cal N_Z(L) ={\Cal Hom}_{\Cal O_H-\text{mod}}(\wedge^n q_*\Cal O_Z,
\wedge^n q_*L). $$

Let $\Cal S=\{(s_{\alpha},L_{\alpha})\}_{\alpha\in\Cal A}$ be a
collection of sections of invertible sheaves on a scheme $Z$, and let
$q : Z \to H$ be a morphism of schemes flat, finite and of relative
rank $n$. We call ${\Cal N}_Z({\Cal S})$ the corresponding
collection of norms on $H$ where 
$$\Cal N_Z(\Cal S)=\{ ( \text{det}(s_{\alpha}),\Cal
N_Z(L_{\alpha}))\}_{\alpha \in \Cal A}.$$

\proclaim{4.5. Proposition} 
Let $q: Z \to H$ be a morphism of schemes, flat, finite and of
relative rank $n$. Let $\Cal S$ be a collection of sections of
invertible sheaves on $Z$, and let $\Cal N_Z(\Cal S)$ be the
corresponding collection of norms on $H$. A morphism of schemes $f :
T\to H$ factors via $\Cal N_Z{(\Cal S)}^{-1}H \to H$ if and only if
the induced morphism of schemes $T\times_H{\Cal S}^{-1}Z \to
T\times_HZ$ is an isomorphism.
\endproclaim

\demo{Proof} This is a global version of Proposition (4.2).\qed
\enddemo

\head  5. - An application to Hilbert schemes of points\endhead

We will in this last section apply results from the previous two
sections about the generalized fraction rings to show the existence of
Hilbert scheme of points on localized schemes ${\Cal S}^{-1}X$, with
$X$ quasi-projective. We will use the fact that the Hilbert scheme of
quasi-projective schemes $X$ exists.

\subhead 5.1. 
Set up\endsubhead We fix a morphism of schemes $X\to S$, where we
refer to $S$ as the base scheme. Let $H$ be an $S$-scheme, and let
$Z\subseteq X\times_S H$ be closed subscheme such that the projection
$q : Z \to H$ is flat, finite and of relative rank $n$. Let $p : Z\to
X$ denote the other projection.

If ${\Cal S}$ is a collection of sections and invertible sheaves on
$X$ we get by the construction (4.4) a collection ${\Cal N}={\Cal
N}_Z(p^*{\Cal S})$ on $H$. We thus have the following diagram
$$\CD
@. (p^{*}\Cal S)^{-1}Z @>>> {\Cal S}^{-1}X \cr
@. @V{i_{p^*{\Cal S}}}VV @V{i_{\Cal S}}VV \cr
 Z_{\Cal N}@>>>  Z @>p>>  X \cr
@VVV @VqVV @. \cr
{\Cal N}^{-1}H @>^{i_{\Cal N}}>> H@.\cr
\endCD \tag{5.1.1}
$$
where the upper right square in (5.1.1) is a fiber product by Lemma
(3.3), and where the scheme $Z_{\Cal N}$ is defined as the fiber
product of the diagram to the down left.

\proclaim{5.2. Lemma} 
The scheme $Z_{\Cal N}$ in the diagram (5.1.1) is a closed subscheme
of ${\Cal S}^{-1}X\times_S {\Cal N}^{-1}H$.
\endproclaim

\demo{Proof} 
We have that $Z$ is closed in $X\times _S H$. It follows that ${\Cal
S}^{-1}X\times_XZ$ is closed in ${\Cal S}^{-1}X\times_SH$. We have
that ${\Cal S}^{-1}X\times_XZ=(q^{*}\Cal S)^{-1}Z$, and thence that
$$ (q^{*}\Cal S)^{-1}Z\times_H {\Cal N}^{-1}H \tag{5.2.1}$$
is a closed subscheme of ${\Cal S}^{-1}X\times_S{\Cal N}^{-1}H$. The
morphism $Z_{\Cal N} \to H$ factors via ${\Cal N}^{-1}H$, hence we
obtain by Proposition (4.4) that (5.2.1) is canonically isomorphic to
$Z_{\Cal N}$. We then have that $Z_{\Cal N}$ is a closed subscheme of
${\Cal S}^{-1}X\times_S{\Cal N}^{-1}H$ as claimed. \qed
\enddemo

\proclaim{5.3. Lemma} Let $f : T \to H$ be a morphism of schemes, and
let $Z_T=Z\times_H T$, where $Z$ a closed subscheme of
$X\times_SH$. If $Z_T\subseteq {\Cal S}^{-1}X\times_SH$ then we have
that the natural morphism $(p^*{\Cal S})^{-1}Z\times_HT \to Z_T$ is an
isomorphism.
\endproclaim

\demo{Proof} If $Z_T$ is a subscheme of ${\Cal S}^{-1}X\times_ST$ then
we have that the projection morphism $Z_T \to X$ factors via the
morphism ${\Cal S}^{-1}X \to X$. Hence $Z_T \to Z$ factors via the
fiber product $(p^*\Cal S)^{-1}Z$ and we obtain our isomorphism.
\enddemo

\definition{5.4. Definition} 
The Hilbert functor ${\Cal H}^n_{X}$ of $n$-points on $X$ is defined
(\cite{G} p. 274) as the contravariant functor from the category of
schemes over $S$ to sets, sending a $S$-scheme $T$ to the set
$$ \align {\Cal H}^n_{X}(T) = \{&\text{closed subschemes } Z\subseteq
X\times_ZT \text{ such that the projection} \\
&\text{ map } q: Z \to T \text{ is flat, finite and of relative rank
}n.\}
\endalign $$
\enddefinition

\proclaim{5.5. Theorem} Let $X\to S$ be a fixed scheme, and assume that the
Hilbert functor of $n$-points on $X$ is represented by a scheme
$\text{H}^n_X$ with universal family $ Z \to \text{H}^n_X$. Let $p :
Z\to X$ denote the projection to $X$. For any collection ${\Cal S}$ of
sections of invertible sheaves on $X$, we let ${\Cal N}=\Cal
N_Z(p^*(\Cal S))$ be the corresponding collection of norms on
$\text{H}^n_X$. We have that the scheme ${\Cal N}^{-1}\text{H}^n_X$ is
the Hilbert scheme of $n$-points on ${\Cal S}^{-1}X$. The universal
family $Z_{\Cal N} \to {\Cal N}^{-1}\text{H}^n_X$ is the pull-back of
the family $Z \to \text{H}^n_X$ along the localization map ${\Cal
N}^{-1}\text{H}^n_X \to \text{H}^n_X$.
\endproclaim
\demo{Proof} 
By Lemma (5.2) we have that $Z_{\Cal N}$ is an ${\Cal
N}^{-1}\text{H}^n_X$-valued point of the Hilbert functor of $n$-points
on ${\Cal S}^{-1}X$. We then have a morphism of functors from the
point functor of ${\Cal N}^{-1}\text{H}^n_X$ to the Hilbert functor
${\Cal H}^n_{\Cal S^{-1}X}$ of $n$-points on ${\Cal S}^{-1}X$. A
morphism we claim is an isomorphism.

Let $T\to S$ be a morphism of schemes, and let $W$ be an $T$-valued
point of ${\Cal H}^n_{\Cal S^{-1}X}$. It follows by Proposition (3.10)
that $W$ is a $T$-valued point of ${\Cal H}^n_X$. Hence there exists a
morphism $f : T \to \text{H}^n_X$ such that the pull-back of the
universal family $Z\to
\text{H}^n_X$ along $f$ is the scheme $W$. We will show that $f$
factors via ${\Cal N}^{-1}\text{H}^n_X$. 

As $W$ is a $T$-valued point of ${\Cal H}^n_{\Cal S^{-1}X}$ it is in
particular a closed subscheme of ${\Cal S}^{-1}X \times_S T$. Hence by
Lemma (5.3) we have that $W=(p^*\Cal S)^{-1}Z
\times_{\text{H}^n_X}T$. That is the natural map
$$ 
(p^*\Cal S)^{-1}Z\times_{\text{H}^n_X}T \to Z \times
_{\text{H}^n_X}T\tag{5.5.1}
$$
is an isomorphism. By Proposition (4.5) the isomorphism (5.5.1) is
equivalent with $f : T\to \text{H}^n_X$ factoring via ${\Cal
N}^{-1}\text{H}^n_X\to \text{H}^n_X$. We thus obtain an morphism of
functors from ${\Cal H}^n_{\Cal S^{-1}X}$ to the point functor of
${\Cal N}^{-1}\text{H}^n_X$, a morphism that clearly is an inverse to
the morphism of functors obtained by the ${\Cal
N}^{-1}\text{H}^n_X$-valued point $Z_{\Cal N}$. \qed
\enddemo

\proclaim{5.6. Corollary} 
Let $U_{\alpha} \subseteq X$ be the open subscheme defined by the
non-vanishing of the section $s_{\alpha} : \Cal O_X \to L_{\alpha}$,
for each $\alpha \in \Cal A$. Then the Hilbert scheme of $n$-points on
$\bigcap_{\alpha \in \Cal A} U_{\alpha}$ is the corresponding
intersection
$$ 
\bigcap_{\alpha \in \Cal A} \text{Hilb}^n_{U_{\alpha}},
$$
where $\text{Hilb}^n_{U_{\alpha}}\subseteq \text{Hilb}^n_X$ is the
open subscheme parameterizing $n$-points on $U_{\alpha}\subseteq X$.
\endproclaim

\demo{Proof} We have that ${\Cal S}^{-1}X=\bigcap _{\alpha \in \Cal
A}U_{s_{\alpha}}$. It then follows from the Theorem that the Hilbert
scheme of points on $\bigcap_{\alpha \in \Cal A}\Cal S^{-1}X$ is the
infinite intersection
$$ 
\bigcap _{\alpha \in \Cal A}H_{\det(s_{\alpha})},
$$ 
where $H_{\det(s_{\alpha})}$ is the open subscheme of
$\text{Hilb}^n_X$ defined by the non-vanishing of the section
$\det(s_{\alpha}) : {\Cal O}_{\text{Hilb}^n_X} \to \Cal N_Z(p^*{\Cal
L_{\alpha}})$. Applying the Theorem to the single pair $(s_{\alpha},
L_{\alpha})$ we get that $H_{\det(s_{\alpha})}$ is the Hilbert scheme
$\text{Hilb}^n_{U_{s_{\alpha}}}$ of $n$-points on
$U_{s_{\alpha}}$. \qed
\enddemo

\subsubhead 5.6.1 
Noetherian schemes \endsubsubhead The Hilbert functor defined in (5.4)
restricts to a functor of noetherian schemes over a noetherian base
scheme $S$.

\proclaim{5.7. Corollary} 
Let $X\to S$ be a projective morphism of noetherian schemes. Let $P\in
X$ be a point, and let ${\Cal O}_{X,P}$ denote the stalk of the
point. Then the Hilbert functor of $n$-points on $\Spec ({\Cal
O}_{X,P})$ is represented by a noetherian scheme. Furthermore, if we
let $\{U_{\alpha}\}$ denote the set of open subschemes $U_{\alpha}$ in
$X$, containing the point $P$, then we have that
$$ 
\operatorname{Hilb}^n_{\Cal O_{X,P}}=\bigcap_{\alpha}
\operatorname{Hilb}^n_{U_{\alpha}}.
$$
\endproclaim

\demo{Proof}  
We have \cite{G} that the Hilbert functor of $n$-points on $X$ is
represented by a projective and in particular noetherian, scheme
$\operatorname {Hilb}^n_{X/S}$. As $X$ is projective we can always
find an locally principal open affine subscheme $U\subseteq X$
containing the point $P$. Hence we have that $\operatorname
{Hilb}^n_{U/S}$, the Hilbert scheme of $n$-points on $U$, is an open
subscheme of $\operatorname{Hilb}^n_{X/S}$. We have that the basic
open affines $D(f)$ form a basis for the topology on $U=\Spec (A)$,
hence we can replace
$$ 
U\cap \bigcap _{\alpha \in \Cal A}U_{\alpha}= \bigcap_{f\in A,P\in
D(f)} D(f).
$$ 
Thus $\Spec ({\Cal O}_{X,P})$ is the localized scheme ${\Cal
S}^{-1}U\subseteq U$, where $S$ the collection $\{f,{\Cal O}_U\}$,
with $P\in D(f)$. We then have by the Theorem that the Hilbert scheme
of $n$-points on ${\Cal S}^{-1}U$ is the infinite intersection of
locally principal subschemes $\bigcap
\text{Hilb}^n_{D(f)}$. The only thing we need to verify is that the
the scheme $\bigcap
\text{Hilb}^n_{D(f)}$ is noetherian. This follows from Proposition (2.7).
\qed
\enddemo

\remark{5.7.1. Remark} 
The Hilbert schemes of points on localized schemes are not generally
varieties, even if the base scheme $S=\Spec (k)$ is the spectrum of a
field. The resulting Hilbert schemes are not always finite type over
the base, and consequently the underlying geometry is complicated if
not bizarre (see \cite{LS}).
\endremark

\remark{5.7.2. Remark} 
Note that the point $P\in X$ in the Corollary, is not assumed to be a
closed point. Thus for an integral scheme $X$ the result also
describes the Hilbert scheme of points on $\Spec (K_X)$, where $K_X$
is the function field of $X$.
\endremark

\enddocument

\end